



\magnification=1100
\baselineskip=15pt 


\catcode`@=11 

\font\ninerm=cmr9
\font\eightrm=cmr8
\font\sixrm=cmr6

\font\ninei=cmmi9
\font\eighti=cmmi8
\font\sixi=cmmi6
\skewchar\ninei='177 \skewchar\eighti='177 \skewchar\sixi='177

\font\ninesy=cmsy9
\font\eightsy=cmsy8
\font\sixsy=cmsy6
\skewchar\ninesy='60 \skewchar\eightsy='60 \skewchar\sixsy='60

\font\eightss=cmssq8

\font\ninebf=cmbx9
\font\eightbf=cmbx8
\font\sixbf=cmbx6

\font\ninett=cmtt9
\font\eighttt=cmtt8

\font\ninesl=cmsl9
\font\eightsl=cmsl8

\font\nineit=cmti9
\font\eightit=cmti8


\newskip\ttglue
\def\tenpoint{\def\rm{\fam0\tenrm}%
  \textfont0=\tenrm \scriptfont0=\sevenrm \scriptscriptfont0=\fiverm
  \textfont1=\teni \scriptfont1=\seveni \scriptscriptfont1=\fivei
  \textfont2=\tensy \scriptfont2=\sevensy \scriptscriptfont2=\fivesy
  \textfont3=\tenex \scriptfont3=\tenex \scriptscriptfont3=\tenex
  \def\it{\fam\itfam\tenit}%
  \textfont\itfam=\tenit
  \def\sl{\fam\slfam\tensl}%
  \textfont\slfam=\tensl
  \def\bf{\fam\bffam\tenbf}%
  \textfont\bffam=\tenbf \scriptfont\bffam=\sevenbf
   \scriptscriptfont\bffam=\fivebf
  \def\tt{\fam\ttfam\tentt}%
  \textfont\ttfam=\tentt
  \tt \ttglue=.5em plus.25em minus.15em

  \normalbaselineskip=12pt
  \def\MF{{\manual META}\-{\manual FONT}}%
  \let\sc=\eightrm
  \let\big=\tenbig
  \setbox\strutbox=\hbox{\vrule height8.5pt depth3.5pt width\z@}%
  \normalbaselines\rm}

\def\ninepoint{\def\rm{\fam0\ninerm}%
  \textfont0=\ninerm \scriptfont0=\sixrm \scriptscriptfont0=\fiverm
  \textfont1=\ninei \scriptfont1=\sixi \scriptscriptfont1=\fivei
  \textfont2=\ninesy \scriptfont2=\sixsy \scriptscriptfont2=\fivesy
  \textfont3=\tenex \scriptfont3=\tenex \scriptscriptfont3=\tenex
  \def\it{\fam\itfam\nineit}%
  \textfont\itfam=\nineit
  \def\sl{\fam\slfam\ninesl}%
  \textfont\slfam=\ninesl
  \def\bf{\fam\bffam\ninebf}%
  \textfont\bffam=\ninebf \scriptfont\bffam=\sixbf
   \scriptscriptfont\bffam=\fivebf
  \def\tt{\fam\ttfam\ninett}%
  \textfont\ttfam=\ninett
  \tt \ttglue=.5em plus.25em minus.15em
  \normalbaselineskip=11pt
  \def\MF{{\manual hijk}\-{\manual lmnj}}%
  \let\sc=\sevenrm
  \let\big=\ninebig
  \setbox\strutbox=\hbox{\vrule height8pt depth3pt width\z@}%
  \normalbaselines\rm}

\def\eightpoint{\def\rm{\fam0\eightrm}%
  \textfont0=\eightrm \scriptfont0=\sixrm \scriptscriptfont0=\fiverm
  \textfont1=\eighti \scriptfont1=\sixi \scriptscriptfont1=\fivei
  \textfont2=\eightsy \scriptfont2=\sixsy \scriptscriptfont2=\fivesy
  \textfont3=\tenex \scriptfont3=\tenex \scriptscriptfont3=\tenex
  \def\it{\fam\itfam\eightit}%
  \textfont\itfam=\eightit
  \def\sl{\fam\slfam\eightsl}%
  \textfont\slfam=\eightsl
  \def\bf{\fam\bffam\eightbf}%
  \textfont\bffam=\eightbf \scriptfont\bffam=\sixbf
   \scriptscriptfont\bffam=\fivebf
  \def\tt{\fam\ttfam\eighttt}%
  \textfont\ttfam=\eighttt
  \tt \ttglue=.5em plus.25em minus.15em
  \normalbaselineskip=9pt
  \def\MF{{\manual opqr}\-{\manual stuq}}%
  \let\sc=\sixrm
  \let\big=\eightbig
  \setbox\strutbox=\hbox{\vrule height7pt depth2pt width\z@}%
  \normalbaselines\rm}

\def\tenbig#1{{\hbox{$\left#1\vbox to8.5pt{}\right.\n@space$}}}
\def\ninebig#1{{\hbox{$\textfont0=\tenrm\textfont2=\tensy
  \left#1\vbox to7.25pt{}\right.\n@space$}}}
\def\eightbig#1{{\hbox{$\textfont0=\ninerm\textfont2=\ninesy
  \left#1\vbox to6.5pt{}\right.\n@space$}}}


\def\titleenv{\def\rm{\fam0\frtnb}%
  \textfont0=\frtnb \scriptfont0=\tenbf \scriptscriptfont0=\sevenbf
  \textfont1=\frtnbmi \scriptfont1=\tenbmi \scriptscriptfont1=\sevenbmi
  \textfont2=\frtnbsy \scriptfont2=\tenbsy \scriptscriptfont2=\sevenbsy
  \textfont3=\tenex \scriptfont3=\tenex \scriptscriptfont3=\tenex
  \def\it{\fam\itfam\frtnbit}%
  \textfont\itfam=\frtnbit
  \def\sl{\fam\slfam\frtnbsl}%
  \textfont\slfam=\frtnbsl
  \def\bf{\fam\bffam\frtnrm}
  \textfont\bffam=\frtnrm \scriptfont\bffam=\tenrm
   \scriptscriptfont\bffam=\sevenrm
  \def\tt{\fam\ttfam\frtntt}%
  \textfont\ttfam=\frtntt
  \tt \ttglue=.5em plus.25em minus.15em

  \normalbaselineskip=17pt
  \def\MF{{\manual META}\-{\manual FONT}}%
  \let\sc=\frtnsc
  \let\big=\tenbig
  \setbox\strutbox=\hbox{\vrule height8.5pt depth3.5pt width\z@}%
  \normalbaselines\rm}


\font\frtnb=cmbx12 scaled\magstep1
\font\frtnbmi=cmmib10 scaled\magstep2
\font\frtnbit=cmbxti10 scaled\magstep2
\font\frtnbsl=cmbxsl10 scaled\magstep2
\font\frtnbsy=cmbsy9 scaled\magstep2
\font\frtnrm=cmr12 scaled\magstep1
\font\frtntt=cmtt10 scaled\magstep2
\font\frtnsc=cmcsc10 scaled\magstep2


\font\tenbmi=cmmib10
\font\sevenbmi=cmmib7
\font\tenbsy=cmbsy9
\font\sevenbsy=cmbsy7


\def\ftnote#1{\edef\@sf{\spacefactor\the\spacefactor}#1\@sf
      \insert\footins\bgroup\eightpoint
      \interlinepenalty100 \let\par=\endgraf
        \leftskip=\z@skip \rightskip=\z@skip
        \splittopskip=10pt plus 1pt minus 1pt \floatingpenalty=20000
        \smallskip\item{#1}\bgroup\strut\aftergroup\@foot\let\next}
\skip\footins=12pt plus 2pt minus 4pt 
\dimen\footins=30pc 


\newcount\ftnoteno
\def\Fnote{\advance\ftnoteno by 1 \ftnote{$^{\number\ftnoteno}$}}


\input amssym.def \input amssym.tex         


\font\sc=cmcsc10




.263   


\def\head#1#2{\headline{\eightss {\ifnum\pageno=1 \underbar{\raise2pt
 \line {#1 \hfill #2}}\else\hfil \fi}}}

\def\title#1{\centerline{\titleenv #1}}

\def\author#1{\bigskip\centerline{\sc By #1}}

\def\abstract#1{\vskip.6in\begingroup\ninepoint\narrower\narrower
 \noindent{\bf Abstract.} #1\bigskip\endgroup}

\def\bottom#1#2#3{{\eightpoint\parindent=0pt\parskip=2pt\footnote{}
 {{\it 1991 Mathematics Subject Classification.}
 #1 }\footnote{}{{\it Key words and phrases.} #2 }\footnote{}{#3 }}}

\outer\def\proclaim #1. #2\par{\medbreak \noindent {\sc #1.\enspace }
 \begingroup\it #2 \endgroup
 \par \ifdim \lastskip <\medskipamount \removelastskip \penalty
 55\medskip \fi}  

\def\section#1. #2.{\vskip0pt plus.3\vsize \penalty -150 \vskip0pt
 plus-.3\vsize \bigskip\bigskip \vskip \parskip
 \centerline {\bf \S#1. #2.}\nobreak \medskip \noindent}



\def\proof{\medbreak\noindent{\it Proof.\enspace}}
\def\qed{\quad \vrule height7.5pt width4.17pt depth0pt}
\def\Qed{\qed\ifmmode \relax \else \medbreak \fi}
\def\remark{\noindent{\sc Remark.\enspace}}
\def\Remark#1\par{\medbreak\remark #1 \medbreak}
\def\example{\noindent{\sc Example:\enspace}}
\def\Example#1\par{\medbreak\example #1 \medbreak}

\def\references{\vskip0pt plus.3\vsize \penalty -150 \vskip0pt
 plus-.3\vsize \bigskip\bigskip \vskip \parskip
 \begingroup\baselineskip=12pt\frenchspacing
   \centerline{\bf REFERENCES}
   \vskip12pt}
\def\ref#1 (#2){\par\hangindent=.8cm\hangafter=1\noindent {\sc #1}\ (#2).}
\def\rfr#1 (#2){\par\hangindent=.8cm\hangafter=1\noindent {\sc #1}\ (#2).}
\def\and{{\rm and }}
\def\endreferences{\bigskip\endgroup}


\def\st{\, ; \;}  
\def\cbuldot{{\raise.25ex\hbox{$\scriptscriptstyle\bullet$}}}


\def\la#1{\mathop{#1}\limits^{\lower 6pt\hbox{$\leftarrow$}}}
\def\ra#1{\mathop{#1}\limits^{\lower 6pt\hbox{$\rightarrow$}}}


\overfullrule=0pt

\input Ref.macros

\theoremcountingtrue
\forwardreferencetrue
\checkdefinedreferencetrue
\sectionnumberstrue
\initialeqmacro

\input psfig

\def\I#1{{\bf 1}_{#1}}
\def\Qed{\qed\ifmmode \relax \else \medbreak \fi}

\font\addrfont=cmcsc10 at 10 truept

\def\E{{\bf E}}

\def\P{{\bf P}}
\def\Q{{\bf Q}}
\def\F{{\cal F}}    
\def\f{{\bf f}}    
\def\G{{\cal Y}}    
\def\Larrow#1{{\buildrel{\lower1.5pt\hbox{$\scriptscriptstyle\leftarrow$}}
 \over {#1}}}
\def\Rarrow#1{{\buildrel{\lower1.5pt\hbox{$\scriptscriptstyle\rightarrow$}}
 \over {#1}}}
\def\longRarrow#1{{\buildrel{\lower1.5pt\hbox{$\scriptscriptstyle
  \longrightarrow$}} \over {#1}}}
\def\Example#1\par{\medbreak\example #1 \medbreak}
\def\Remark#1\par{\medbreak\remark #1 \medbreak}

\def\Lh{\widehat L}
\def\Xh{\widehat X}
\def\Th{\widehat T}
\def\GW{{\bf GW}}  
\def\GWh{{\widehat{{\bf GW}}}}  

\def\nuh{\widehat \nu}

\def\ugwh{\GWh{}_*}  

\ifproofmode \relax \else\head{To appear in Ann. Probab.}{Version of 19 September 1994}\fi

\vglue20pt
\title{Conceptual Proofs of $L \log L$ Criteria}
\title{for Mean Behavior of Branching Processes}
\author{Russell Lyons, Robin Pemantle, and Yuval Peres}
\bottom{Primary 60J80.}{Galton-Watson, 
size-biased distributions.}
{Research partially supported by two Alfred P. Sloan Foundation
Research Fellowships (Lyons and Pemantle), by NSF Grants DMS-9306954
(Lyons), DMS-9300191 (Pemantle),
and DMS-9213595 (Peres), by a Presidential Faculty
 Fellowship (Pemantle),
 and by the hospitality of Washington University (Pemantle).}

\abstract{The Kesten-Stigum Theorem is a fundamental criterion for the rate
of growth of a supercritical branching process, showing that an $L \log L$ 
condition is decisive. In critical and subcritical cases, 
 results of Kolmogorov and later authors
 give the rate of decay of the probability that the
process survives at least $n$ generations.  
We  give conceptual proofs of these theorems based on comparisons
of Galton-Watson measure to another measure on the space of trees.
This approach also explains Yaglom's
exponential limit law for conditioned critical branching processes via
a simple characterization of the exponential distribution.}

\nobreak
\bsection{Introduction}{s.intro}

Consider a Galton-Watson branching process with each particle having
probability $p_k$ of generating $k$ children. Let $L$ stand for a random
variable with this offspring distribution. Let  $m := \sum_k k p_k$ be the
mean number of children per particle and let $Z_n$ be the number 
of particles in the $n^{th}$ generation.  The most basic and well-known fact
about branching processes is that the extinction probability $q := 
\lim \P [Z_n = 0]$ is equal to 1 if and only if $m \le 1$ and $p_1 < 1$.  
It is also not hard to establish that in the case $m > 1$,
$$
{1 \over n}\log Z_n  \to \log m
$$ 
almost surely on nonextinction, while in the case $m \le 1$, 
$$
{1 \over n}\log \P [Z_n > 0]  \to \log m .
$$
 Finer questions may be asked:
\smallskip
\item{$\bullet$} In the case $m > 1$, when does the mean
$\E[Z_n]= m^n$ give the right growth rate up to a random factor?  
\item{$\bullet$} In the case $m < 1$, when
does the first moment estimate $\P [Z_n > 0] \le \E[Z_n]= m^n$ give
the right decay rate up to a random factor?
\item{$\bullet$} In the case $m = 1$, 
what is the decay rate of $\P [Z_n > 0]$? 
\smallskip\noindent
These questions
are answered by the following three classical theorems.

\proclaim Theorem A: Supercritical Processes (Kesten and Stigum (1966)).
 \hfill \break
Suppose that
$1 < m < \infty$ and let $W$ be the limit of the martingale $Z_n/m^n$.
The following are equivalent:
\smallskip
\item{(i)}  $\P[W=0] = q$;
\item{(ii)} $\E[W] = 1$;
\item{(iii)}  $\E[L\log^+ L] < \infty$.

\medskip


\proclaim Theorem B: Subcritical Processes 
(Heathcote, Seneta and Vere-Jones (1967)). \hfil\break
The sequence $\{ \P [ Z_n > 0]/m^n \}$ is decreasing. If $m < 1$, then
the following are equivalent:
\item{(i)} $\lim_{n\to\infty} \P [ Z_n > 0]/m^n > 0$;
\item{(ii)} $\sup \E[Z_n \mid Z_n > 0] < \infty$;
\item{(iii)} $\E [L \log^+ L] < \infty$.


The fact that (i) holds if $\E[L^2] < \infty$
 was proved  by Kolmogorov (1938).
 It is interesting that 
the law of $Z_n$ conditioned on $Z_n > 0$
always converges in a strong sense, even
when its means are unbounded; 
see \ref s.williamson/.

\medskip

\proclaim  Theorem C: Critical Processes (Kesten, Ney and Spitzer (1966)). 
\hfill \break  Suppose that $m = 1$ and let $\sigma^2 := {\rm Var}(L) =
\E[L^2] - 1 \le \infty$.  Then we have
\smallskip
\item{(i)} Kolmogorov's estimate:
$$\lim_{n \to \infty} n \P [ Z_n > 0] = {2 \over \sigma^2}\,;$$
\item{(ii)} Yaglom's limit law: 
\item{} If $\sigma < \infty$, then
 the conditional distribution of ${Z_n/n}$ given $Z_n>0$ converges 
 as $n \rightarrow \infty$ to an exponential law with mean ${\sigma^2}/2 \,$.
If $\sigma = \infty$, then this conditional distribution converges to
infinity.
\medskip

Under a third moment assumption, parts  (i) and (ii) of Theorem C are due to
Kolmogorov (1938) and Yaglom (1947), respectively.

 For classical proofs of these theorems, the reader is referred to Athreya
and Ney (1972), pp.~15--33 and 38--45 or Asmussen and Hering (1983),
pp.~23--25, 58--63, and 74--76. A very short proof of the
Kesten-Stigum theorem, using martingale truncation, is in Tanny (1988).

By using simple measure theory, we reduce the
dichotomies between mean and sub-mean behavior
 in the first two theorems to easier known dichotomies concerning the growth
of branching processes with immigration. These, in turn, arise from
the following dichotomy, which is an immediate consequence of the 
Borel-Cantelli lemmas.

\procl l.edich Let $X, X_1, X_2, \ldots$ be nonnegative i.i.d.\ random
variables.
Then
$$
\limsup_{n \to \infty} {1 \over n} X_n = \cases{0 &if \ $\E[X] < \infty\, ,$\cr
    \infty &if \ $\E[X] = \infty \, .$ \cr}
$$
\endprocl
\smallskip

Size-biased distributions, which arise in many contexts, 
play an important role in the present paper.
 Let $X$ be a nonnegative
random variable with finite positive mean. Say that $\Xh$ has the
 corresponding {\it  size-biased} distribution if
$$\E[g(\Xh)]={{\E[Xg(X)]} \over{\E[X]}}$$   
 for every positive Borel
 function $g$. The analogous notion for random trees is the topic of
\ref s.size-biased/. 

Note that if $X$ is an exponential random variable
and $U$ is uniform in $[0,1]$ and independent of $\Xh$, then the 
product $U \cdot \Xh$ has the same distribution as $X$. 
(One way to see this is by considering the first and second points
of a Poisson process.)
 The fact that
this property actually characterizes the exponential distributions (Pakes
and Khattree (1992)) is used
in \ref s.Kolm/ to derive part (ii) of Theorem C.

The next section is  basic for the rest of the paper;
Sections 3, 4 and 5, which contain the proofs of Theorems A, C and B,
respectively, may be read independently of each other.
An extension of Theorem A to branching processes in a random
environment, due to Tanny (1988), is  discussed in the final section.

\bsection{Size-biased Trees}{s.size-biased}

Our proofs depend on viewing Galton-Watson processes as generating random
 family trees, not merely as generating
various numbers of particles; of course, this 
goes back at least to Harris (1963).
  We think of these trees as rooted and labeled,
with the (distinguishable) offspring of each vertex ordered from left to right.
We shall define another way of growing random
trees, called {\bf size-biased Galton-Watson}.
 The law of this random tree
will be denoted $\GWh$, whereas the law of an ordinary Galton-Watson tree
is denoted $\GW$.

Let $\Lh$ be a random variable
whose distribution is that of size-biased $L$,  i.e.,
$\P[\Lh = k] = {k p_k / m}$.
To construct a size-biased Galton-Watson tree $\Th$, start with an 
initial particle
$v_0$.  Give it a random number $\Lh_1$ of children, where $\Lh_1$ has the
law of $\Lh$. Pick one of these children at random, $v_1$. Give the {\it
other\/} children independently
 ordinary Galton-Watson descendant trees and give $v_1$ an independent
size-biased number $\Lh_2$ of children. Again, pick one of the children
of $v_1$ at random, call it $v_2$, and give the others ordinary Galton-Watson
descendant trees. Continue in this way indefinitely.
(See Figure \figref{sizeGW}.) Note that
size-biased Galton-Watson trees are always infinite (there is no
extinction).  


Define the measure $\ugwh$ as the joint
distribution of the random tree $\Th$ and the random path 
$ (v_0,v_1,v_2,\ldots)$. Let $\GWh$ be its marginal on the space of
trees.

 For a tree $t$ with $Z_n$ vertices at level $n$,
write $ W_n(t) := Z_n/m^n$. For any rooted tree $t$ and any $n \ge 0$, denote
by $[t]_n$ the set of rooted trees whose first $n$ levels agree with those
of $t$. (In particular, if the height of $t$ is less than $n$, then $[t]_n
= \{t\}$.) If $v$ is a vertex at the $n$th level of $t$, then
let $[t;v]_n$ denote the set of 
{\bf trees with distinguished paths} such that the tree is in $[t]_n$
and the path starts from the root, does not backtrack, and goes through $v$. 

 Assume that $t$ is a tree
of height at least $n+1$ and that the root of $t$ has $k$ children with
descendant trees $t^{(1)}, t^{(2)}, \ldots , t^{(k)}$.
Any vertex $v$ in level $n+1$ of $t$ is in one of these,
say $t^{(i)}$. 
The measure $\ugwh$ clearly satisfies the recursion
$$
\ugwh[t;v]_{n+1} = {{k p_k} \over {m}} \cdot {1 \over k} \cdot
             \ugwh [t^{(i)};v]_{n} \cdot 
                    \prod_{j \neq i}\GW[t^{(j)}]_{n} \, .
$$ 
By induction, we conclude that
$$
 \ugwh [t;v]_n ={1 \over {m^n}} \GW[t]_n   \label e.goal
$$
 for all $n$ and all $[t;v]_n$ as above.
Therefore,
$$
\GWh[t]_n = W_n(t) \GW[t]_n \, , \label e.ac
$$
 for all $n$ and all trees $t$.
 From \ref e.goal/ we see
that, given the first $n$ levels of the tree $\Th$, the measure 
$\ugwh$ makes the  vertex $v_n$
in the random path $(v_0, v_1,\ldots)$  uniformly distributed
on the $n$th level of $\Th$.

The vertices off the ``spine'' $(v_0,v_1,\ldots)$ of the size-biased tree
 form a {\bf branching process with immigration}. 
In general, such a process is defined by two distributions, an offspring 
distribution and an immigration distribution. The process
starts with no particles, say, and at every generation $n \ge 1$, there is an 
immigration of $Y_n$ particles, where $Y_n$ are i.i.d.\ with the given
immigration law.
Meanwhile, each particle has, independently, an ordinary Galton-Watson 
descendant tree with the given offspring  distribution.

 Thus, the $\GWh$-law of $Z_n - 1$ is the 
same as that of
the generation sizes of an immigration process with $Y_n = \Lh_{n} - 1$.
 The probabilistic content of the assumption $\E [L \log^+ L]
< \infty$ will arise in applying \ref l.edich/ to the
variables $\{ \log^+ Y_n \}$,
since $\E[ \log^+ (\Lh-1)] = m^{-1} \E [L \log^+ (L-1)]$.

\bigskip
The construction of size-biased trees is not new. It and related
constructions in other situations occur in Kahane and Peyri\`ere (1976),
Kallenberg (1977),
Hawkes (1981), Rouault (1981),
Joffe and Waugh (1982), Kesten (1986),
Chauvin and Rouault (1988), Chauvin, Rouault and
Wakolbinger (1991), and Waymire and Williams (1993). The paper of Waymire
and Williams (1993) is the only one among these to use such a
construction in a similar way to the method we use to prove Theorem A;
their work was independent of and contemporaneous with ours.
None of these papers use methods similar to the ones we employ for the
proofs of Theorems B and C.
 An {\it a priori\/} motivation for the use of size-biased trees in our
context comes from 
the general principle that in order to study asymptotics,
it is useful to construct a suitable limiting object first.
In the supercritical case, to study the asymptotic behavior
of the martingale $W_n := Z_n / m^n$ with  respect to $\GW$, it is natural
to consider the  sequence of measures $\, W_n \,d \GW$, which 
converge weakly to $\GWh$. As pointed out by the editor, this can also be
viewed as a Doob $h$-transform. When $m \leq 1 $, the size-biased tree
may be obtained by conditioning a Galton-Watson tree
to survive forever. 
 The generation sizes of
size-biased Galton-Watson trees are known as a Q-process 
in the case $m \le 1$; see Athreya-Ney (1972), pp.\ 56-60.
One may also view $\ugwh$ as a Campbell measure and $\GWh$ as the
associated Palm measure.

\bsection{Supercritical Processes: Proof of Theorem A}{s.KS}

Theorem A will be an immediate consequence of the following theorem on 
immigration processes.

\procl t.seneta \procname{Seneta (1970)} Let $Z_n$ be the generation sizes
of a Galton-Watson process with immigration $Y_n$ . Let $m := \E[L]>1$ be the
mean of the offspring  law and let $Y$ have the same 
law as $Y_n$. If $\E[\log^+ Y] < \infty$, then $\lim Z_n /m^n$ exists and is finite
a.s., while if $\E[\log^+ Y] = \infty$, then $\limsup Z_n/c^n = \infty$ a.s. for 
every constant $c>0$.
\endprocl


\proof (Asmussen and Hering (1983), pp.~50--51) Assume
 first that $\E[\log^+ Y] = \infty$. By \ref l.edich/, $\limsup
Y_n/c^n = \infty$ a.s. Since $Z_n \ge Y_n$, the result follows.

Now assume that $\E[\log^+ Y] < \infty$.
Let $\G$ be the $\sigma$-field generated by $\{Y_k \st k \ge 1
\}$.   Let $Z_{n,k}$ be the number of
descendants at level $n$ of the vertices which immigrated in generation $k$.
 Thus, the total number of vertices at level $n$ is
$\sum_{k=1}^{n} Z_{n,k}$. This gives
$$
 \E[Z_n/m^n \mid \G] =  \E\left[ {1 \over m^n} \sum_{k=1}^{n} Z_{n,k}
   \,\bigg | \,\G \right]
= \sum_{k=1}^{n} {1 \over m^k} \E \left[ { Z_{n,k} \over  m^{n-k}}
   \,\bigg |\, \G \right] \, . 
$$
Now, for $k < n$, the random variable 
 $Z_{n,k} / m^{n-k}$ is the $(n-k)$th element of the
ordinary Galton-Watson martingale sequence starting with,
however, $Y_k$ particles. Therefore, its expectation is just
$Y_k$ and so
$$
\E[Z_n/m^n \mid \G]  = \sum_{k=1}^{n} {Y_k \over m^k}  \, .
$$
Our assumption gives, by \ref l.edich/, that $Y_k$ grows
subexponentially, whence this series converges a.s.
Since $\{Z_n/m^n\}$ is a submartingale when
conditioned on $\G$ with bounded expectation (given $\G$), it converges
a.s. \Qed

To prove Theorem A,
recall the following elementary result, whose proof we include for the sake
of completeness:

\procl p.Wdich  Either $W=0$ a.s.\ or $W>0$ a.s.\ on nonextinction.
In other words, ${\P[W =0]\in\{q,1\}}$.
\endprocl

\proof Let $f(s) := \E[s^L]$ be the probability generating function of $L$.
The roots of $f(s) = s$ in ${[0,\, 1]}$ are $\{q, 1\}$. Thus, it suffices to
show
that $\P[W=0]$ is such a root. Now the $i$th individual of the first
generation has a descendant Galton-Watson tree with, therefore, a
martingale limit, $W^{(i)}$, say. These are independent and have the same
distribution as $W$. Furthermore,
$$
W = {1\over m} \sum^{Z_1}_{i=1} W^{(i)}, 
$$
or, what counts for our purposes,
$$
W =0\iff \forall i\le Z_1\ \  W^{(i)} =0\, .
$$
Conditioning on $Z_1$ now gives immediately the desired fact that
$f(\P[W=0]) = \P[W=0]$. 
\Qed

\noindent{\it Proof of Theorem A.}\enspace
Rewrite \ref e.ac/ as follows.
Let $\F_n$ be the $\sigma$-field generated by the first $n$ levels
of trees and $\GW_n$, $\GWh_n$ be the restrictions of $\GW$, $\GWh$ to
$\F_n$. Then \ref e.ac/ is the same as
$$
{d\GWh_n \over d\GW_n} (t) = W_n(t) \, .\label e.rnn
$$
It is convenient now to interpret the last expression
 for {\it infinite} trees $t$, where both sides
depend only on the first $n$ levels of $t$.
In order to define $W$ for every infinite tree $t$, set
$$
W(t) := \limsup_{n\to\infty} W_n(t) \, .
$$
 From \ref e.rnn/ follows the key dichotomy:
$$
W=0 \quad \GW\hbox{-a.s.} \iff \GW \perp \GWh \iff W
= \infty\quad \GWh\hbox{-a.s.} \label e.di1
$$
 while 
$$
\int W \,d\GW = 1 \iff \GWh \ll \GW \iff W
< \infty  \quad \GWh\hbox{-a.s.} \label e.di2
$$
 (see Durrett (1991), p.\ 210, Exercise 3.6).
This is the key because it
allows us to change the problem from one about the $\GW$-behavior of
$W$ to one about the $\GWh$-behavior of $W$.
Indeed, since the $\GWh$-behavior of $W$ is described by \ref t.seneta/,
the theorem is immediate: if $\E[L \log^+ L] <  \infty$, i.e.,  $\E[\log^+ \Lh] < \infty$,
then $W < \infty $ $\GWh$-a.s. by \ref t.seneta/,
whence $\int W \,d\GW = 1$ by \ref e.di2/;
while if $\E[L \log^+ L] = \infty$, then $W = \infty$ $\GWh$-a.s. by \ref
t.seneta/, whence $W = 0$ $\GW$-a.s. by \ref e.di1/. \Qed
\vskip -\parskip
\vskip -\baselineskip
\bsection{Critical Processes: Proof of Theorem C}{s.Kolm}



\procl l.con Consider a critical Galton-Watson process with 
a random number $Y\ge 1$ of initial particles in generation 0.
 Choose one of the initial
particles, $v$, at random. Let $B_n$ be the
event that at least one of the particles to the left of $v$ has a
descendant in generation $n$ and let
$r_n$ be the number of descendants in
generation $n$ of the particles to the right of $v$.
Let $\beta_n := \P[B_n]$ and $\alpha_n := \E[r_n \I{ B_n}]$. Then 
$\lim_{n \to \infty} \beta_n  = 0$
and, if $\E[Y] < \infty$, then
$\lim_{n \to \infty} \alpha_n = 0$.
\endprocl

\proof The fact that $\beta_n \to 0$ follows from writing $\P[B_n] =
\E\big[\P[B_n \mid Y, v]\big]$ and applying the bounded convergence theorem.
Now
$$
\alpha_n = \E\big[\E[r_n \I{B_n} \mid Y, v]\big]
 = \E\big[\E[r_n \mid Y, v] \P[B_n \mid Y, v] \big]
 \le \E\big[Y \P[B_n \mid Y, v] \big] = \E[Y \I{B_n}]
$$
by independence of $r_n$ and $B_n$ given $Y$ and $v$. Thus, $\E[Y] <
\infty$ implies that $\alpha_n \to 0$. \Qed


\noindent{\it Proof of Theorem C (i)}. \enspace
Let $A_n$ be the event that $v_n$ is the leftmost vertex in generation $n$.
By definition, $\ugwh (A_n \mid Z_n) = 1 / Z_n$.  From this, it follows
that conditioning on $A_n$ reverses the effect of size-biasing.  
That is, the law of the first $n$ generations of a tree under $(\ugwh \mid
A_n)$ is the same as under $(\GW \mid Z_n > 0)$. In particular,
$$\int Z_n \, d (\ugwh \mid A_n) = \int Z_n \, d(\GW \mid Z_n > 0) 
   = {1 \over \GW (Z_n > 0)} .$$

We are thus required to show that
$$
{1 \over n}\int Z_n \, d (\ugwh \mid A_n) \to {\sigma^2 \over 2}\, .\label
e.x
$$
 For any tree with a distinguished line of descendants $v_0 , v_1 , 
\ldots$, decompose the size of the $n$th generation by writing
$Z_n = 1 + \sum_{j=1}^{n} Z_{n,j}$, 
where $Z_{n,j}$ is the number of vertices at generation $n$
descended from $v_{j-1}$ but not from $v_{j}$.
The intuition
behind \ref e.x/ is that the unconditional $\ugwh$-expectation of $Z_{n,j}$
is $\E[\Lh]- 1 = \sigma^2$; half of these fall to the left of
$v_n$ and half to the right. Since the chance that any given vertex at
generation $n-k$ other than $v_{n-k}$ has no
 descendant in generation $n$ tends to 1 as
$k \to \infty$,  conditioning on none surviving
to the left leaves us with $\sigma^2 / 2$. 

To prove this, define $R_{n, j}$ to be the number of vertices in
generation $n$ descended from those children of $v_{j-1}$ to the right of
$v_{j}$ and $R_n := 1 + \sum_{j=1}^{n} R_{n, j}$, the number of vertices
in generation $n$ to the right of $v_n$, {\it inclusive}.
Let $A_{n,j}$ be the event that
$R_{n,j} = Z_{n,j}$.  
Let $R'_{n,j}$ be independent random variables
with respect to a probability measure $\Q'$ such that $R'_{n,j}$ has
the $(\ugwh \mid A_{n,
j})$-distribution of $R_{n, j}$. Let $\Q := \ugwh \times \Q'$.  Define 
$$
R^*_{n, j} := R_{n, j} \I{A_{n ,j}} + R'_{n,j} \I{\neg A_{n,j}} \,,
$$
where $\neg$ denotes complement.
Then the random variable $R^*_n := 1 + \sum_{j=1}^{n}
R^*_{n, j}$ has the same distribution as the $(\ugwh \mid A_n)$-law of
$Z_n$ since the event $A_n=\, \{R_n=Z_n\}$ is the
 intersection of the independent events $A_{n,j}$. Also, 
$$
\int R_{n, j} \,d(\ugwh \mid A_{n, j})
  \le \int Z_{n, j} \,d(\ugwh \mid A_{n, j})
  \le \int Z_{n,j} \,d\ugwh = \E[\Lh]-1 = \sigma^2\,,
$$
where the second inequality is due to
$Z_{n,j}$ and the indicator of $A_{n,j}$ being negatively correlated.
Now, for each $j$,
we  apply \ref l.con/ with $Y = \Lh_j$ to the descendant
trees of the children of $v_{j-1}$, with $\neg A_{n,j}$ playing the role
of $B_{n-j}$. We conclude that if $\sigma < \infty$,
then
\begineqalno
\int {1 \over n} |R_n - R^*_n| \,d\Q &\le \int {1 \over n} \sum_{j=1}^{n}
|R_{n,j} - R^*_{n,j}|\,d\Q \cr
 & \label e.same \cr
 &\le  {1 \over n} \sum_{j=1}^{n} 
    \int_{\neg A_{n,j}} (R_{n,j} + R^*_{n,j}) \,d\Q
 \le {1 \over n} \sum_{j=1}^{n} (\alpha_{n-j}  + \sigma^2 \beta_{n-j})
 \to 0 \cr 
\endeqalno
as $n \to \infty$.
In particular, since $\int R_{n,j} \,d\ugwh = \sigma^2/2$
and hence $\int R_n/n \,d\ugwh = \sigma^2/2$,
we get $\int R^*_n/n\,d\Q  \to \sigma^2/2$. The case 
$\sigma = \infty$ follows from this
by truncating $\Lh_k$ while leaving unchanged the rest of the
size-biased tree. This shows \ref e.x/, as desired. \Qed

The following simple characterization of the exponential distributions
is used to prove part (ii) of the theorem.

\procl l.char \procname{Pakes and Khattree (1992)}
Let $X$ be a nonnegative random variable with a positive finite mean
and let $\Xh$ have the corresponding size-biased distribution.
Denote by $U$ a uniform random variable in $[0,1]$ which is independent
of $\Xh$.
Then $U \cdot \Xh$ and $X$ have the same distribution iff $X$ is exponential.
\endprocl

\proof  By linearity, we may assume that $\E[X] = 1$.
 For any $\lambda > 0$, we have
$$
\E\left[e^{-\lambda U \cdot \Xh} \right] 
= \E \left[ \int_0^1 X e^{-\lambda u \cdot X} \,du \right]
= {1 \over \lambda} \E\left[ 1 - e^{-\lambda X} \right] \,,
$$
which equals $\E[e^{-\lambda X}]$ iff $\E[e^{-\lambda X}] = 1/(\lambda +
1)$. By uniqueness of the Laplace transform, this holds for every $\lambda > 0$
iff $X$ is exponential with mean 1. \Qed

The following lemma is elementary.

\procl l.sz Suppose that $X$, $X_n$ are nonnegative random variables
with positive finite means such that $X_n
\to X$ in law and $\Xh_n \to Y$ in law. 
If $Y$ is a proper random variable, then $Y$ has the law of $\Xh$.
\endprocl

\medskip
\noindent{\it Proof of Theorem C (ii)}. \enspace
Suppose first that $\sigma < \infty$. This ensures that the $\GWh$-laws of
$Z_n/n$ have uniformly
bounded means and, {\it a fortiori}, are tight. Let $R_n$ and
$R^*_n$ be as in the proof of part (i).
Then $R^*_n/n$ also have uniformly bounded means and hence are tight.
Therefore, there is a sequence $\{n_k\}$ tending
to infinity such that $R^*_{n_k}/n_k$ converges in law to a (proper)
random variable $X$ and the $\GWh$-laws of $Z_{n_k}/n_k$ converge to
the law of a (proper) random variable $Y$. 
Note that the law of $R^*_{n_k}$ is the $(\GW \mid Z_n > 0)$-law of
$Z_{n_k}$. Thus, from
\ref l.sz/ combined with \ref e.x/, the variables $Y$ and $\Xh$ are identically
distributed. Also, by \ref e.same/,
the $\ugwh$-laws of $R_{n_k}/n_k$ tend to the law of $X$.

  On the other
hand, let $U$ be a uniform [0, 1]-valued random variable independent of
every other random variable encountered so far.  Then $R_n$ and $\lceil
U\cdot Z_n \rceil$ have the same law (with respect to $\ugwh$), while
$$
\left| {1 \over n} \lceil U\cdot Z_n \rceil
   - {1 \over n} U\cdot Z_n \right| \le {1 \over n} \to 0\,.
$$
Hence $X$ and $U\cdot \Xh$ have the same distribution.
It follows from \ref l.char/ and \ref e.x/ that $X$ is
an exponential random variable with mean $\sigma^2/2$. In particular, the
limiting distribution of $R^*_{n_k}/n_k$ is
independent of the sequence $n_k$, and hence we actually have convergence
in law of the whole sequence $R^*_n/n$ to $X$, as desired.

Now suppose that $\sigma = \infty$. A truncation argument shows that the
$\GWh$-laws of $Z_n/n$ tend to infinity, whence so do the laws of $\lceil
U \cdot Z_n \rceil / n$. Thus, the $(\GW \mid Z_n > 0)$-laws of $Z_n/n$ tend
to infinity as well. \Qed

\Remark The fact that the limit $\GWh$-law of $Z_n/n$ is that of $\Xh$, i.e.,
the sum of two independent exponentials with mean $\sigma^2/2$ each, is due
to Harris (see Athreya and Ney (1972), pp.~59--60). The above proof allows
us to identify these two exponentials as normalized counts
of the vertices to the left and right of the ``spine'' $(v_0, v_1, \ldots)$.

\bsection{Subcritical Processes: Proof of Theorem B}{s.Yaglom}

Let $\mu_n$ be the law of $Z_n$ conditioned on $Z_n > 0$.
 For any tree $t$, let $\xi_n (t)$ be the leftmost vertex in the first 
generation having at least one descendant in generation $n$ if $Z_n > 0$.
Let $H_n (t)$ be the number of descendants of $\xi_n (t)$ in generation 
$n$, or zero if $Z_n = 0$.  It is easy to see that
$$
\GW(H_n = k \mid Z_n > 0) = \GW (H_n = k \mid Z_n > 0 ,\, \xi_n = x) = 
   \GW (Z_{n-1} = k \mid  Z_{n-1} > 0)
$$
 for all children $x$ of the root.  Since $H_n \le Z_n$, this shows that 
$\{\mu_n\}$ increases stochastically in $n$.
Now
$$
\GW ( Z_n > 0) = {\E[Z_n] \over \E [Z_n \mid Z_n > 0]} = {m^n \over \int
x\,d\mu_n(x)} \,.
$$
Therefore, $\GW(Z_n>0)/m^n$ is decreasing and (i) $\Leftrightarrow$ (ii).
The equivalence of (ii) and (iii) is an immediate consequence of
the following routine lemma applied to the laws $\mu_n$ and the
 following theorem on immigration. \Qed

\procl l.tightness Let $\{ \nu_n \}$ be a sequence of probability measures
on the positive integers with finite means $a_n$.  Let $\nuh_n$ be
size-biased, i.e., $\nuh_n (k) = k \nu_n (k) / a_n$.  If
   $\{ \nuh_n \}$ is tight,
then $\sup a_n < \infty$, while if
$\nuh_n \to \infty$ in distribution, then $a_n \to \infty$.
\endprocl


\procl t.heath \procname{Heathcote (1966)} Let $Z_n$ be the generation sizes
of a Galton-Watson process with offspring random variable $L$ and
immigration $Y_n$ . Suppose that
$m := \E[L]<1$  and let $Y$ have the same 
law as $Y_n$.  If $\E[\log^+ Y] < \infty$, then 
$Z_n$ converges in distribution to a proper random variable,
while if $\E[\log^+ Y] = \infty$, then $Z_n$ converges in probability to infinity.
\endprocl

The following proof is a slight improvement on Asmussen and Hering (1983),
pp.~52--53.

\proof Let
$\G$ be the $\sigma$-field generated by $\{Y_k \st k \ge 1
\}$.  For any $n$, let $Z_{n,k}$ be the number of
descendants at level $n$ of the vertices which immigrated in generation $k$.
 Thus, the total number of vertices at level $n$ is
$\sum_{k=1}^{n} Z_{n,k}$. Since the distribution of $Z_{n,k}$ depends only on
$n-k$, this total $Z_n$ has the same distribution as 
$ \sum_{k=1}^n Z_{2k,k} $,
which is an increasing process with limit $Z'_\infty$. By Kolmogorov's zero-one
law, $Z'_\infty$ is a.s.\ finite or a.s.\ infinite.
Hence, we need only to show that $Z'_\infty < \infty$ iff $\E[\log^+ Y] <
\infty$.

Assume that $\E [\log^+ Y] < \infty$.  
Now $\E[Z'_\infty \mid \G] = \sum_{k=1}^\infty Y_k m^{k-1}$.
Since $\{ Y_k \}$ is almost surely subexponential in $k$ by \ref l.edich/,
this sum converges a.s.
Therefore,  $Z'_\infty$ is finite a.s.

Now assume that $Z'_\infty < \infty$ a.s. Writing
$Z_{2k,k} = \sum_{i=1}^{Y_k}  \zeta_k(i)$, where $ \zeta_k(i)$ are
the sizes of generation $k-1$ of
i.i.d.\ Galton-Watson branching processes with one initial particle,
we have $Z'_\infty = \sum_{k=1}^\infty
\sum_{i=1}^{Y_k}  \zeta_k(i)$ written as a random
sum of independent random
variables. Only a finite number of them are at least one, whence
by the Borel-Cantelli lemma conditioned on $\G$, we get
$\sum_{k=1}^\infty Y_k \GW(Z_{k-1} \ge 1) < \infty$ a.s.
Since $\GW(Z_{k-1} \ge 1) \ge \P[L > 0]^{k-1}$, it follows 
by \ref l.edich/ that $\E[\log^+ Y] < \infty$. \Qed
\vskip -\parskip
\vskip -\baselineskip
\bsection{Strong Convergence of the Conditioned Process in the Subcritical 
          Case}{s.williamson}

Yaglom (1947) showed that when $m<1$ and $Z_1$ has a finite second moment,
the conditional distribution $\mu_n$ of $Z_n$ given $\{Z_n>0\}$
converges to a proper probability distribution as $n \rightarrow \infty$.
This was proved without the second moment assumption by Joffe (1967)
and by Heathcote, Seneta and Vere-Jones (1967). The following stronger
convergence result was proved, in an equivalent form, by Williamson
(cf. Athreya-Ney (1972), pp. 64--65.)

\procl t.will The sequence $\{\mu_n\}$ always converges in a strong sense:
if $\| \cdot \|$ denotes total variation norm, then
$\sum_n \|\mu_n - \mu_{n-1}\| < \infty$.
\endprocl

\Remark Note that this is strictly stronger than weak
convergence to a probability measure, even for a stochastically increasing
sequence of distributions.

\proofof t.will
Recalling the notation of the previous section and the
events $A_{n,j}$ from \ref s.Kolm/, we see that
$$
{1 \over 2} \|\mu_n - \mu_{n-1}\| \le \GW (H_n \neq Z_n \mid Z_n > 0) 
  = \ugwh(H_n \neq Z_n
  \mid A_n) = \ugwh (H_n \neq Z_n \mid A_{n,1}) \,.
$$
Let $\lambda$ be the number of children of the root to the left 
of $v_1$ and let $s_n = \GW (Z_n > 0)$.  
Now condition on $\Lh_1$ and $\lambda$ and use the fact that 
$\inf_{n \ge 2} \ugwh (A_{n,1}) =: \delta > 0$ to estimate
$$\eqalign{
\ugwh (H_n \neq Z_n \mid A_{n,1})
   &\le \delta^{-1} \ugwh (\{ H_n \neq Z_n \} \cap A_{n,1}) \cr
   &= \delta^{-1} \sum_{k=1}^\infty \sum_{l=0}^{k-1} \ugwh (\Lh_1 = k ,
      \lambda = l , H_n \neq Z_n, A_{n,1}) \cr
   &= \delta^{-1} \sum_{k=1}^\infty {k p_k \over m} \sum_{l=0}^{k-1} {1 \over k}
       (1 - s_{n-1})^l [1 - (1 - s_{n-1})^{k-1-l}] \,.
}$$
Sum this in $n$ by breaking it into two pieces: those $n$ for which
$s_{n-1}^{-1} \le k$ and those for which $s_{n-1}^{-1} > k$.  For 
the first piece, use 
$\sum_{l=0}^{k-1} (1 - s_{n-1})^l \le s_{n-1}^{-1}$,
and for the second piece, use 
$$
\sum_{l=0}^{k-1} [1 - (1 - s_{n-1})^{k-1-l}]
\le \sum_{l=0}^{k-1} (k-1-l) s_{n-1} \le k^2 s_{n-1} / 2\,.
$$ 
These estimates yield
$$\sum_{n=1}^\infty \|\mu_n - \mu_{n-1}\|
   \le 2\delta^{-1} \sum_k {p_k \over m} \left[
      \sum_{s_{n-1}^{-1} \le k} s_{n-1}^{-1} + 
      \sum_{s_{n-1} < 1/k} k^2 s_{n-1}/2 \right] .$$
By virtue of Theorem B, we have $s_j \le m s_{j-1}$, 
so that each of these two inner sums is bounded by
 a geometric series, whence the total is finite. \Qed

\bsection{Stationary Random Environments}{s.notes}

The analogue of the Kesten-Stigum theorem for branching processes in random
environments (BPRE's) is due to Tanny (1988). 
Our method of proof applies
to this situation too.
 Furthermore, the probabilistic construction of size-biased trees
makes apparent how to remove
 a technical hypothesis in Tanny's extension of the
Kesten-Stigum theorem. 
The analytic tool 
needed for this turns out to be an ergodic lemma due to Tanny (1974).

In this set-up, the fixed ``environment'' $f$ is replaced by a stationary
ergodic sequence $f_n$ of random environments. Vertices in generation $n-1$
have offspring  according to the law with p.g.f.\ $f_n$. Assume that the
process is supercritical with finite growth, i.e., $0 < \E[\log f'_1(1)] <
\infty$.  Write $M_n:=
\prod_{k=1}^{n} f'_k(1)$; this is the conditional mean of $Z_n$ given the
environment sequence $\f := \langle f_k \rangle$. The quotients $Z_n/M_n$ still
 form a martingale with a.s.\ limit $W$. 
 We shall also use our notation
that $\Lh_n$ are random variables which, given $\f$, are independent and have
size-biased distributions, so that $f'_n(s) = f'_n(1) \sum_{k \ge 1}
\P[\Lh_n = k] s^{k-1}$.
Tanny's (1988) theorem with the
technical hypothesis removed is as follows.

\procl t.tanny If for some $a>0$, the sum $\sum_n \P[\Lh_n > a^n \mid \f\,] $
is finite a.s., then $\E[W]=1$ and $W \ne 0$ a.s.\ on nonextinction, while
if this sum is infinite with positive probability
 for some $a > 0$, then $W = 0$ a.s. In case the
environments $f_n$ are i.i.d., the a.s.\ finiteness of this sum for some $a$ is
equivalent to the finiteness of $\E[\log^+ \Lh_1]$.
\endprocl

In the proof of \ref t.tanny/, for the case of i.i.d.\ environments,
\ref l.edich/ is used in the same way as it is
 for a fixed environment. However, the general
case requires the following extension of \ref l.edich/.

\procl l.tannyerg \procname{Tanny (1974)} Let
$\tau$ be an ergodic measure-preserving
transformation of a probability space $(X,\, \mu)$, and let $f$ be a nonnegative
measurable function on  $X$. Then 
$\limsup f(\tau^n(x))/n$ is either 0 a.s.\ or $\infty$ a.s.
\endprocl

The following argument, indicated to us by Jack Feldman, is considerably
shorter than the proofs in Tanny (1974) and  Wo\'s (1987). It is similar
to but somewhat shorter than the proof in O'Brien (1982).

\proofof l.tannyerg 
 By ergodicity, $\limsup f(\tau^n(x))/n$ is a.s.\ a constant $c \leq \infty$.
 Suppose that
 $0 <c< \infty$. Then we may choose $L$ so that
 $A:=\{x \st \forall k \ge L \  f(\tau^k(x))/k <2c\}$
 has $\mu(A) >{9 \over 10}$. The ergodic theorem applied
 to $\I{A}$ implies
 that for a.e.\ $x$, if $n$ is sufficiently large, then there exists 
  $k \in [L,\, n/5]$ such that $y:=\tau^{n-k}(x) \in A$, and therefore
$$
{f(\tau^n(x)) \over n} = {f \left(\tau^k(y) \right) \over k}
  \cdot {k \over n} < {2c \over 5} \, .
$$
This contradicts the definition of $c$. \Qed

The role of \ref l.tannyerg/ in the proof
of \ref t.tanny/ is not large:  By \ref l.tannyerg/, we know that
$\limsup (1/n) \log^+
\Lh_n$ is a.s.\ 0 or a.s.\ infinite. Since the random variables
$\Lh_n$ are independent given the environment $\f$, the Borel-Cantelli
lemma still shows that which of these alternatives holds
is determined by whether the sum $\sum_n \P[\Lh_n > a^n \mid \f\,] $ is
 finite for every $a > 0$ or not. In particular,
this sum is finite with positive probability
 for some $a$ if and only if it is finite a.s.\ for every $a$.

The proof that if $\limsup (1/n) \log^+ \Lh_n = \infty$ a.s., then $W=0$
a.s.\ applies without change to the case of random environments. The proof
in the other direction needs only conditioning on $\f$ in
addition to conditioning on the $\sigma$-field $\G$.

The fact that if $\E[W] > 0$, then $W\ne 0$ a.s.\ on nonextinction for
BPRE's can be proved
by virtually the same method as that of \ref p.Wdich/, using
the uniqueness of the conditional extinction probability (Athreya and
Karlin 1971).

\bigskip

\bigskip
\noindent{\bf Acknowledgement:} We are grateful to Jack Feldman for permitting
the inclusion of his proof of \ref l.tannyerg/ in this paper.


\references

\rfr Asmussen, S. \and Hering, H. (1983) Branching Processes. Birkh\"auser,
Boston.

\rfr Athreya, K. B. \and Karlin, S. (1971) On branching processes with
random environments I: Extinction probabilities, {\it Ann. Math. Statist.}
{\bf 42}, 1499--1520.

\rfr Athreya, K. B. \and Ney, P. (1972) Branching Processes. Springer, New
York.


\rfr Chauvin, B. \and Rouault, A. (1988) KPP equation and supercritical
branching Brownian motion in the subcritical speed area. Application to
spatial trees, {\it Probab. Theory Relat. Fields} {\bf 80}, 299--314.

\rfr Chauvin, B., Rouault, A., \and Wakolbinger, A. (1991) Growing
conditioned trees, {\it Stochastic Process. Appl.} {\bf 39}, 117--130.

\rfr Durrett, R. (1991) Probability: Theory and Examples. Wadsworth,
Pacific Grove, California.

\rfr Harris, T. E. (1963) The Theory of Branching Processes. Springer, Berlin.

\rfr Hawkes, J. (1981) Trees generated by a simple branching process, {\it
J. London Math. Soc.} {\bf 24}, 373--384.

\rfr Heathcote, C. R. (1966) Corrections and comments on the paper ``A
branching process allowing immigration", {\it J. Royal Statist. Soc. {\bf
B}} {\bf 28}, 213--217.

\rfr Heathcote, C. R., Seneta, E., \and Vere-Jones, D. (1967)
A refinement of two theorems in the theory of branching processes,
{\it Theory Probab. Appl.} {\bf 12}, 297--301.

\rfr Joffe, A. (1967) On the Galton-Watson branching processes with
mean less than one, {\it Ann. Math. Statist.} {\bf 38}, 264--266.

\rfr Joffe, A. \and Waugh, W. A. O'N. (1982) Exact distributions of kin
numbers in a Galton-Watson process, {\it J. Appl. Prob.} {\bf 19},
767--775.

\rfr Kahane, J.-P. \and Peyri\`ere, J. (1976) Sur certaines martingales
de Benoit Mandelbrot, {\it Adv. in Math.} {\bf 22}, 131--145.

\rfr Kallenberg, O. (1977) Stability of critical cluster fields, {\it Math.
Nachr.} {\bf 77}, 7--43.

\rfr Kesten, H. (1986) Subdiffusive behavior of random walk on
a random cluster, {\it Ann. Inst. H. Poincar\'e Probab. Statist.}
 {\bf 22}, 425--487.

\rfr Kesten, H., Ney, P. \and Spitzer, F. (1966) The Galton-Watson
process with mean one and finite variance, {\it Theory Probab. Appl.}
{\bf 11}, 513--540.

\rfr Kesten, H. \and Stigum, B. (1966) A limit theorem for multidimensional
Galton-Watson processes, {\it Ann. Math. Statist.} {\bf 37}, 1211--1223.

\rfr Kolmogorov, A. N. (1938) On the solution of a problem in biology,
{\it Izv. NII Matem. Mekh. Tomskogo Univ.} {\bf 2}, 7--12.



\rfr O'Brien, G. L. (1982) The occurrence of large values in stationary
sequences, {\it Z. Wahrsch. Verw. Gebiete} {\bf 61}, 347--353.

\rfr Pakes, A. G. \and Khattree, R. (1992) Length-biasing,
characterization of laws and the moment problem, {\it Austral. J.
Statist.} {\bf 34}, 307--322.

\rfr Rouault, A. (1981) Lois empiriques dans les processus de branchement
spatiaux homog\`enes supercritiques, {\it C. R. Acad. Sci. Paris S\'er.
I. Math.} {\bf 292}, 933--936.

\rfr Seneta, E. (1970) On the supercritical branching process with immigration,
 {\it Math. Biosci. } {\bf 7}, 9--14.


\rfr Tanny, D. (1974) A zero-one law for stationary sequences, {\it
Z. Wahrsch. Verw. Gebiete} {\bf 30}, 139--148.

\rfr Tanny, D. (1988) A necessary and sufficient condition for a branching
process in a random environment to grow like the product of its means, {\it
Stochastic Process. Appl.} {\bf 28}, 123--139.

\rfr Waymire, E. C. \and Williams, S. C. (1993) A cascade decomposition
theory with applications to Markov and exchangeable cascades, {\it preprint}.

\rfr Wo\'s, J. (1987) The filling scheme and the ergodic theorems of
Kesten and Tanny, {\it Colloq. Math.} {\bf 52}, 263--276.

\rfr Yaglom, A. M. (1947) Certain limit theorems of the theory
of branching processes, {\it Dokl. Acad. Nauk. SSSR} {\bf 56}, 795--798.

\endreferences

\begingroup
\addrfont
\parindent=0pt\baselineskip=10pt

Department of Mathematics,
Indiana University,
Bloomington, IN 47405-5701

Department of Mathematics,
University of Wisconsin,
Madison, WI 53706

Department of Statistics,
University of California,
Berkeley, CA 94720

\endgroup

\bye